\documentclass[10pt,a4paper,reqno]{amsart}

\usepackage{comment,graphicx,amssymb,xcolor,amsmath,mathtools,enumerate} 
\usepackage{hyperref}
\theoremstyle{plain}
\newtheorem{theorem}{Theorem}%[section]
\newtheorem{lemma}[theorem]{Lemma}
\newtheorem{proposition}[theorem]{Proposition}
\newtheorem{question}[theorem]{Question}
        
\theoremstyle{remark}
\newtheorem{example}[theorem]{Example}
\newtheorem{remark}[theorem]{Remark}
        
\numberwithin{equation}{section}

\begin{document}

\title[Irrationality of rapidly converging series]{Irrationality of rapidly converging series: \\ A problem of Erd\H{o}s and Graham}

% author information
\author[K.~Barreto]{Kevin Barreto}
\address[Kevin Barreto]{Queens' College, University of Cambridge, Cambridge, CB3 9ET, United Kingdom}
\email{kb799@cam.ac.uk}
%\urladdr{}

% author information
\author[J.~Kang]{Jiwon Kang}
\address[Jiwon Kang]{Department of Mathematical Science, Seoul National University, Seoul, 08826, Korea}
\email{carl011226@gmail.com}
%\urladdr{}

% author information
\author[S.~Kim]{Sang-hyun Kim}
\address[Sang-hyun Kim]{School of Mathematics, Korea Institute for Advanced Study (KIAS), Seoul, 02455, Korea}
\email{skim.math@gmail.com}
\urladdr{https://kimsh.kr}

% author information
\author[V.~Kova\v{c}]{Vjekoslav Kova\v{c}}
\address[Vjekoslav Kova\v{c}]{University of Zagreb Faculty of Science, Department of Mathematics, Bijeni\v{c}ka cesta 30, 10000 Zagreb, Croatia}
\email{vjekovac@math.hr}
\urladdr{https://web.math.pmf.unizg.hr/~vjekovac/}

% author information
\author[S.~Zhang]{Shengtong Zhang}
\address[Shengtong Zhang]{Department of Mathematics, Stanford University, CA 94305, USA}
\email{stzh1555@stanford.edu}

%\date{July 2, 2026}


\begin{abstract}
Answering a question of Erd\H{o}s and Graham, we show that the double-exponential growth condition $\limsup_{n\to\infty}a_n^{1/\phi^n}=\infty$ for a strictly increasing sequence of positive integers $\{a_n\}_{n=1}^\infty$ is sufficient for the series $\sum_{n=1}^\infty 1/(a_n a_{n+1})$ to have an irrational sum; here $\phi$ denotes the golden ratio. We also provide a positive generalization to $\sum_{n=1}^\infty 1/(a_n^{w_0}\cdots a_{n+d-1}^{w_{d-1}})$, and a negative result showing that some of these instances are essentially optimal. The original problem was autonomously solved by the AI agent \emph{Aletheia}, powered by Gemini Deep Think, while the remaining material is largely a product of human-AI interactions.
\end{abstract}

\maketitle

\section{Declaration of AI usage}
The present paper resolves a problem posed by Paul Erd\H{o}s and Ronald Graham~\cite{ErdosGraham1980}, which subsequently appeared in~\cite{Erd1988,BloomWeb}.
Beyond the mathematical result, this work is also notable for its methodology, which relies heavily on artificial intelligence in a collaborative effort with Google DeepMind.
We outline the research timeline below:

\begin{enumerate}
\item Erd\H{o}s and Graham posed Question~\ref{que:erdos1051} below in their well-known 1980 problem book~\cite[p.~64]{ErdosGraham1980}, expressing the belief that it should be solvable (perhaps by standard tools) and that some generalization should exist.
\item 
\emph{Aletheia}~\cite{aletheia}\footnote{Aletheia is a custom mathematics research agent built upon Gemini Deep Think, details are available at \url{https://github.com/google-deepmind/superhuman/tree/main/aletheia}.} found a solution to Question~\ref{que:erdos1051}.
From the original output, it was clear that some generalization was possible. Specifically, the proof used rough estimates such as $1.7 < 2$, $\log x < x$, etc. See~\cite{aletheia-github} for relevant raw outputs.
\item The authors realized and proved that one can relax the original growth condition \eqref{eq:Erd_hyp} and generalize the series \eqref{eq:Erd_two_terms}; this generalization is now formulated as Part \eqref{it:T1} of Theorem~\ref{thm:erdos1051-sharp}.
\item Gemini Deep Think discovered a further generalization (among many plausible variations of the hypotheses) and then, jointly with \emph{Aletheia}, proved it. This generalization corresponds to Theorem~\ref{thm:erdos1051-gen}, with $b_n=1$ and a stronger additional hypothesis \eqref{eq:stronger_than_needed}.
\item The authors found counterexamples that tightly complement the previously mentioned result; these are discussed in Part \eqref{it:T2} of Theorem~\ref{thm:erdos1051-sharp}.
\item Finally, inspired by Erd\H{o}s's paper~\cite{Erdos1975}, the authors were able to weaken the hypotheses further to obtain Theorem~\ref{thm:erdos1051-gen} below. In addition, a more general negative result was found and formulated as Theorem~\ref{thm:erdos1051-construction}.
\end{enumerate}

Thus, the theorems and proofs in this paper are largely products of human-AI interactions; each of the steps above also depends on the ideas discovered in the previous steps. 
In particular, the AI discovered the formulation of Theorem~\ref{thm:erdos1051-gen} after the human authors presented Theorem~\ref{thm:erdos1051-sharp}.
This work also represents an interesting use case for AI in that it originated when \emph{Aletheia} searched over a large set of open problems and autonomously narrowed them down to those it could solve.
However, we emphasize that the writing of the paper was entirely done by the (human) authors.

%%%%%

\section{Introduction}
The following question was formulated by Erd\H{o}s and Graham in~\cite[p.~64]{ErdosGraham1980} and posed again by Erd\H{o}s several years later~\cite[p.~106]{Erd1988}. We quote its self-contained formulation from the website \emph{Erd\H{o}s problems}~\cite{BloomWeb}.

\begin{question}\label{que:erdos1051}
Is it true that if $a_1<a_2<\cdots$ is a sequence of positive integers with
\begin{equation}\label{eq:Erd_hyp} 
\liminf_{n\to\infty} a_n^{1/2^n}>1, 
\end{equation}
then 
\begin{equation}\label{eq:Erd_two_terms}
\sum_{n=1}^\infty \frac{1}{a_n a_{n+1}}
\end{equation}
is necessarily irrational?
\end{question}

The present paper provides a positive answer to Question~\ref{que:erdos1051} and its natural generalizations, namely Theorems~\ref{thm:erdos1051-sharp} and~\ref{thm:erdos1051-gen}.
The proposers of Question~\ref{que:erdos1051}, Erd\H{o}s and Graham, noted that hypothesis \eqref{eq:Erd_hyp} ``should be enough'' to imply the conclusion, while remarking that ``it is not too hard to show that if $a_n\to\infty$ rapidly enough''~\cite[p.~64]{ErdosGraham1980}. Their latter remark can be inferred from~\cite[Theorem 1]{Erdos1975} to mean, for example, that $\lim_{n\to\infty} a_n^{1/2^n}=\infty$.
They also asked ``what the strongest theorem of this type'' would be, and our result, Theorem~\ref{thm:erdos1051-sharp} below, provides a satisfactory answer.

Series of the form 
\[ \sum_{n=1}^{\infty} \frac{1}{a_n} \]
for a strictly increasing sequence of positive integers $\{a_n\}_{n=1}^\infty$ are sometimes called \emph{Ahmes series}. The topic of relating their irrationality solely to the growth of the sequence was systematically studied in~\cite{ES1964,Erdos1975}, with very special instances appearing much earlier; see~\cite{KT2025} for some subsequent literature.
On the other hand, series of the form 
\begin{equation}\label{eq:Cantor}
\sum_{n=1}^{\infty} \frac{b_n}{a_1 a_2 \cdots a_n}
\end{equation}
for some positive integer sequences $\{a_n\}_{n=1}^\infty$ and $\{b_n\}_{n=1}^\infty$ with $a_n\geq2$ for every $n$ are known as \emph{Cantor series}. 
Investigation of their irrationality dates back to the mid-19\textsuperscript{th} century~\cite{Cantor1869,Oppenheim1954,ES1971,ES1974}, while more recent relevant literature includes~\cite{TY2002,HT2004,HT2004b}.
These results can be used to solve particular cases of Question~\ref{que:erdos1051} simply by rewriting \eqref{eq:Erd_two_terms} as the Cantor series \eqref{eq:Cantor} with $b_n=a_1\cdots a_{n-2}$.
In this way, a result by Han\v{c}l and Tijdeman~\cite[Theorem 3.2]{HT2004} confirms that \eqref{eq:Erd_two_terms} is not a rational number under the following assumptions on divisibility and the ``relative growth'' among the terms $a_n$:
\begin{align*} 
& a_n \nmid a_1\cdots a_{n-2} \text{ for every } n\in\mathbb{N}, \quad \liminf_{n\to\infty} \frac{a_1\cdots a_{n-2}}{a_n} = 0, \\
& \text{and}\quad \lim_{n\to\infty} \frac{a_1\cdots a_{n-2}}{a_{n-1} a_n} = 0. 
\end{align*}

Our result below confirms and generalizes the irrationality of \eqref{eq:Erd_two_terms} with no assumptions on the relative growth.

\begin{theorem}\label{thm:erdos1051-sharp}
Fix a positive integer $d$ and let $\psi>1$ be the unique positive solution to $\psi^d=\psi^{d-1}+1$.
\begin{enumerate}
\item\label{it:T1}
If a monotonically increasing\footnote{i.e., non-decreasing} sequence of positive integers $\{a_n\}_{n=1}^\infty$ satisfies
\[ \lim_{n \to \infty} a_n^{1/\psi^n} =\infty, \]
then
\begin{equation}\label{eq:w_are_1}
\sum_{n=1}^\infty \frac{1}{a_n a_{n+1} \cdots a_{n+d-1}}
\end{equation}
is an irrational number.
\item\label{it:T2} 
For every $C\in(1,\infty)$ there exists a strictly increasing sequence of positive integers $\{a_n\}_{n=1}^\infty$ satisfying
\begin{equation*}
\lim_{n\to\infty} a_n^{1/\psi^n} = C
\end{equation*}
such that the infinite sum \eqref{eq:w_are_1} is a rational number.
\end{enumerate}
\end{theorem}

When we specialize Theorem~\ref{thm:erdos1051-sharp} by taking $d=2$, $\psi$ becomes the golden ratio $\phi=(1+\sqrt{5})/2<2$. In particular, it implies that condition \eqref{eq:Erd_hyp} is sufficient for the irrationality of \eqref{eq:Erd_two_terms}, whereas the condition $\liminf_{n\to\infty} a_n^{1/\phi^n}>1$ is not sufficient.

The following heuristic can explain the appearance of the equation $\psi^d=\psi^{d-1}+1$ determining the critical rate of growth of the sequence $\{a_n\}_{n=1}^\infty$ that guarantees irrationality of the sum \eqref{eq:w_are_1}. It was sketched by Terence Tao in the case of the golden ratio $\phi$ \cite{BloomWeb}. Assume that, ideally, $a_n$ exhibits a steady double-exponential growth: $a_n \approx \exp(c^n)$ for some number $c>1$. The $N$-th partial sum of \eqref{eq:w_are_1} is a rational number with denominator at most
\[ a_1 a_2 \cdots a_{N+d-1} \approx \exp(c^1+c^2+\cdots+c^{N+d-1}), \]
while the corresponding series tail has magnitude
\[ \approx \frac{1}{\exp(c^{N+1}+c^{N+2}+\cdots+c^{N+d})}, \]
which is determined by the $(N+1)$-st term, since the terms with indices $n\geq N+2$ are relatively negligible. If the whole sum \eqref{eq:w_are_1} happens to be rational, then the tail needs to catch up with the error coming from the $N$-term approximation:
\[ c^{N+1}+c^{N+2}+\cdots+c^{N+d} \leq c^1+c^2+\cdots+c^{N+d-1}. \]
This can be rewritten as
\[ \frac{c^{N+1}(c^d-1)}{c-1} \leq \frac{c^{N+d}-c}{c-1}, \]
easily giving $c^d-1\leq c^{d-1}$, i.e., $c\leq \psi$, where $\psi$ was defined in Theorem~\ref{thm:erdos1051-sharp}. Conversely, if the sequence grows like $a_n \approx \exp(\psi^n)$ or slower, then we can hope that dominant tails will allow us to perturb $\{a_n\}_{n=1}^\infty$ slightly in a way to guide the partial sums towards a nearby rational number. Thus, the growth $\exp(\psi^n)$ really is critical. The actual proof of the theorem will necessarily be more complicated, since the rate of growth of $\{a_n\}_{n=1}^\infty$ can vary.

To the best of our knowledge Theorem~\ref{thm:erdos1051-sharp} is a novel result when $d\geq2$. For $d=1$, we get $\psi=2$, so Part \eqref{it:T1} is an easy consequence of Erd\H{o}s's theorem~\cite[Theorem~1]{Erdos1975}, while the well-known Sylvester sequence provides an explicit example of a sequence as in Part \eqref{it:T2}; see the remarks in~\cite[Section~1]{KT2025}. 
For $d\geq2$, Part \eqref{it:T1} is a special case of an even more general Theorem~\ref{thm:erdos1051-gen} below, while Part \eqref{it:T2} is a particular instance of Theorem~\ref{thm:erdos1051-construction}.

\begin{theorem}\label{thm:erdos1051-gen}
Fix a positive integer $d$ and a tuple of non-negative integers $\mathbf{w} = (w_0, w_1, \dots, w_{d-1})$ such that $w_{d-1} \ge 1$. Set $W = \max\{w_0, \dots, w_{d-1}\}$, and consider the unique positive real root $c_{\mathbf{w}}$ of the polynomial $$ P_{\mathbf{w}}(x) \coloneq  (x-1)\sum_{j=0}^{d-1} w_j x^j - W x^{d-1}. $$ Let $\{a_n\}_{n=1}^\infty$ and $\{b_n\}_{n=1}^\infty$ be sequences of positive integers, with $\{a_n\}_{n=1}^\infty$ monotonically increasing. Assume that there exist real numbers $0<\eta<\tau$ such that 
\begin{equation*} 
b_n\le n^\eta,\quad a_{n}^{w_0}a_{n+1}^{w_1}\cdots a_{n+d-1}^{w_{d-1}}\geq n^{1+\tau}
\end{equation*}
for all $n\in\mathbb{N}$. We further assume that
$$
\limsup_{n\to\infty} a_n^{1/c_{\mathbf{w}}^n} =\infty.$$
Then the infinite sum
    $$ S_{\mathbf{w}}(\{a_n\}_{n=1}^\infty,\{b_n\}_{n=1}^\infty) \coloneq  \sum_{n=1}^\infty \frac{b_n}{a_{n}^{w_0}a_{n+1}^{w_1}\cdots a_{n+d-1}^{w_{d-1}}} $$
is irrational.
\end{theorem}

\begin{remark}
\begin{enumerate}
\item It is elementary to see that $P_{\mathbf{w}}(x)=0$ has a unique positive solution $c_{\mathbf{w}}$, which belongs to $(1,\infty)$. 
Namely, $P_{\mathbf{w}}$ is clearly negative on $(0,1]$, while writing
\[ P_{\mathbf{w}}(x) = (x-1) Q(x) - W, \quad Q(x) \coloneq w_{d-1} x^{d-1} - \sum_{j=0}^{d-2} (W-w_j) x^j \]
and applying Descartes' rule of signs to $Q$, we easily conclude that $P_{\mathbf{w}}$ has precisely one root in $(1,\infty)$.
\item Observe that for every sequence of positive integers $\{a_n\}_{n=1}^\infty$ and
numbers $\vartheta>\theta>1$ we have the implication
\[ \liminf_{n \to \infty} a_n^{1/\vartheta^n} > 1 \quad\implies\quad \lim_{n \to \infty} a_n^{1/\theta^n} = \infty. \]
Hence, the (stronger than needed) condition 
\begin{equation}\label{eq:stronger_than_needed}
\lim_{n\to\infty} a_n^{1/c^n_{\mathbf{w}}}=\infty
\end{equation}
is met whenever we have some $c>c_{\mathbf{w}}$ satisfying
$$\liminf_{n\to\infty} a_n^{1/c^n}>1.$$
\item Another consequence of Theorem~\ref{thm:erdos1051-gen} is the following. Suppose $d\ge2$ and $\psi^d=\psi^{d-1}+1$ for some $\psi>1$. If $\{a_n\}_{n=1}^\infty$ is a strictly increasing sequence of positive integers satisfying $$\limsup_{n\to\infty} a_n^{1/\psi^n}=\infty,$$ then the sum $$\sum_{n=1}^\infty \frac1{a_na_{n+1}\cdots a_{n+d-1}}$$ is irrational. To see this, note first that $a_n\ge n$, and hence,  $${a_n a_{n+1}\cdots a_{n+d-1}}\ge n^d\ge n^2.$$ Thus, the hypotheses of Theorem~\ref{thm:erdos1051-gen} are automatically satisfied. In particular, this proves the result mentioned in the abstract.
\item Erd\H{o}s~\cite[Theorem 1]{Erdos1975} proved the irrationality of $\sum_{n=1}^{\infty} 1/a_n$ under the assumptions
$$ a_n \geq  n^{1+\tau} $$
for some $\tau>0$ and every $n\in\mathbb{N}$ and
$$ \limsup_{n\to\infty} a_n^{1/2^n}=\infty. $$
This is implied by the case $d=1$ in Theorem~\ref{thm:erdos1051-gen}.
Numerous variants and generalizations of this theorem by Erd\H{o}s exist (see, e.g.,~\cite{HS03,HN17,HL19}), but we could not find any concerning series of the form \eqref{eq:Erd_two_terms} or any of the series addressed by cases $d\geq 2$ of Theorem~\ref{thm:erdos1051-gen}.
However, one should note that the aforementioned result of Erd\H{o}s trivially implies the instances of Theorem~\ref{thm:erdos1051-gen} when $c_{\mathbf{w}}\geq 2$, which happens if and only if $\sum_{j=0}^{d-1}2^j w_j\leq 2^{d-1} W$.
\end{enumerate}
\end{remark}

The following theorem provides a negative counterpart to Theorem~\ref{thm:erdos1051-gen}.

\begin{theorem}\label{thm:erdos1051-construction}
Fix a positive integer $d$ and a tuple of non-negative integers $\mathbf{w} = (w_0, w_1, \dots, w_{d-1})$ such that $w_{d-1} \ge 1$. 
If $\tilde{c}_{\mathbf{w}}$ is the largest positive root of the polynomial 
$$ \widetilde{P}_{\mathbf{w}}(x) \coloneq  (x-1)\sum_{j=0}^{d-1} w_j x^j - x^{d-1}, $$ 
then, for every $C\in(1,\infty)$, there exists a strictly increasing sequence of positive integers $\{a_n\}_{n=1}^\infty$ satisfying
\begin{equation}\label{eq:examplegrowth}
\lim_{n\to\infty} a_n^{1/\tilde{c}_{\mathbf{w}}^n} = C
\end{equation}
and
\[ \sum_{n=1}^\infty \frac{1}{a_{n}^{w_0}a_{n+1}^{w_1}\cdots a_{n+d-1}^{w_{d-1}}} \in \mathbb{Q}. \]
\end{theorem}

\begin{remark}
\begin{enumerate}
\item The largest real root of $\widetilde{P}_{\mathbf{w}}$ lies in $(1,\infty)$ because of $\widetilde{P}_{\mathbf{w}}(1)<0$.
\item Note that Theorem~\ref{thm:erdos1051-construction} implies that Theorem~\ref{thm:erdos1051-gen} is sharp in all instances when $w_0, \dots, w_{d-1}\in\{0,1\}$.
\end{enumerate}
\end{remark}

\begin{example}
When $\mathbf{w}= (1,0,2,1)$, the unique positive root of the polynomial
\[ P_{\mathbf{w}}(x) = (x-1)(1+2x^2+x^3)-2x^3 \]
is
\[ c_{\mathbf{w}}=1.914\ldots, \]
while the unique positive root of the polynomial
\[ \widetilde{P}_{\mathbf{w}}(x) = (x-1)(1+2x^2+x^3)-x^3 \]
is
\[ \tilde{c}_{\mathbf{w}}=1.345\ldots. \]
Whenever a monotonically increasing sequence of integers $\{a_n\}_{n=1}^\infty$ satisfies $a_n\geq n^{1+\tau}$ for some $\tau>0$
and $\limsup_{n\to\infty} a_n^{1/c_{\mathbf{w}}^n}=\infty$,
we see from Theorem~\ref{thm:erdos1051-gen} that the sum
$$ \sum_{n=1}^\infty \frac{1}{a_n a_{n+2}^2 a_{n+3}} $$
is irrational. On the other hand, Theorem~\ref{thm:erdos1051-construction} implies that this sum can be a rational number for a strictly increasing sequence of integers $\{a_n\}_{n=1}^\infty$ such that $\lim_{n\to\infty} a_n^{1/\tilde{c}_{\mathbf{w}}^n}\in(1,\infty)$.
\end{example}

%%%%%

\section{Proof of Theorem~\ref{thm:erdos1051-gen}}
We utilize a classical irrationality criterion that dates back to at least Fourier's proof of the irrationality of the number $e$, first published in~\cite{Sta1815}.
This completely elementary lemma is often attributed to Mahler~\cite{Mahler1976,Nishioka1996}.

\begin{lemma}\label{lem:mahler}
Let $\{z_n\}_{n=1}^\infty$ be non-negative rational numbers such that 
\[ S\coloneq\sum_{n=1}^\infty z_n<\infty \] 
and $z_n>0$ for infinitely many $n\in\mathbb{N}$. Let $\{D_N\}_{N=1}^\infty$ be positive integers such that the value $D_N\sum_{n=1}^N z_n$ is an integer for every $N\in\mathbb{N}$. If $S$ is rational, then 
\[ \liminf_{N\to\infty} D_N\Bigl(S-\sum_{n=1}^N z_n\Bigr)>0. \]
\end{lemma}
\begin{proof}
Write $S=p/q$ in lowest terms. Then $D_N\bigl(S-\sum_{n=1}^N z_n\bigr)$ is a positive rational with denominator dividing $q$, hence it is at least $1/q$.
\end{proof}

We also use a simple trick, which Erd\H{o}s attributed to Borel.

\begin{lemma}[{cf.~\cite{Erdos1975}}]\label{lem:borel} 
If $\{\delta_n\}_{n=1}^\infty$ and $\{\mu_n\}_{n=1}^\infty$ are positive sequences with $\sum_{n=1}^\infty \delta_n<\infty$ and $\limsup_{n\to\infty} \mu_n=\infty$,
then the set
\begin{equation*}
\mathcal{P}\coloneq\bigl\{m\in\mathbb{N} \,:\, \mu_{m+1}>(1+\delta_{m})\max_{n\le m}\mu_{n}\bigr\}
\end{equation*}
is infinite.
\end{lemma}

\begin{proof}
Assuming the contrary and denoting the largest element of $\mathcal{P}$ by $n_0\in\mathbb{N}$, we get that every integer $n\ge n_0+1$ satisfies 
\[ \mu_{n+1} \leq \mu_{n_0+1}\prod_{k=n_0+1}^n (1+\delta_{k}) \leq \mu_{n_0+1}\prod_{k=1}^\infty (1+\delta_{k}) .\] 
Since $\prod_{k=1}^{\infty} (1+\delta_k) \leq\exp(\sum_{k=1}^{\infty} \delta_k)<\infty$, we conclude that $\{\mu_n\}_{n=1}^\infty$ is bounded, which contradicts our assumptions.
\end{proof}

The following elementary lemma appears in~\cite{Erdos1975}. We prove it by splitting positive integers into dyadic blocks. Given two functions $f$ and $g$, let us write $f(x)\ll g(x)$ when there exists some constant $C>0$ satisfying $|f(x)|\leq C|g(x)|$ for all $x$ in their common domain. 
This is the Vinogradov notation for $f=O(g)$.

\begin{lemma}\label{lem:integral-weighted-pattern}
Let $e\le x_1\le x_2\le\cdots$ be a monotonically increasing sequence of real numbers and
assume $0\le y_k\le k^\eta$ for some $\eta\ge 0$ and all $k$. Let $1\le Q<R$ be integers.
\begin{enumerate}
\item\label{it:Wpoly-pattern}
If $x_k\ge k^{1+\tau}$ for all $k\in[Q,R]$ with some $\tau>\eta$, then  
$$\sum_{k=Q}^R \frac{y_k}{x_k}\ll \frac{1}{x_Q^{(\tau-\eta)/(1+\tau)}}.$$
\item\label{it:Wexp-pattern}
If $x_k\ge e^k$ for all $k\in[Q,R]$, then 
$$\sum_{k=Q}^R \frac{y_k}{x_k}\ll \frac{(\log x_Q)^{1+\eta}}{x_Q}.$$
\end{enumerate}
The implicit constants here only depend on $\eta$ and $\tau$.
\end{lemma}

\begin{proof}
For every integer $j\ge 0$ set 
$$I(j)\coloneq\{k\in[Q,R]:\ 2^j x_Q\le x_k\le 2^{j+1}x_Q\}.$$ 
Then, 
\begin{equation}\label{eq:auxbx}
\sum_{k=Q}^R\frac{y_k}{x_k}\le \sum_{j=0}^\infty \sum_{k\in I(j)}\frac{y_k}{x_k}\le \sum_{j=0}^\infty \frac{\sum_{k\in I(j)} y_k}{2^j x_Q}.
\end{equation}

First, we prove Part \eqref{it:Wpoly-pattern}.
If $k\in I(j)$, then 
$$y_k^{1+\tau}\le k^{(1+\tau)\eta}\le x_k^\eta\le \left(2^{j+1}x_Q\right)^\eta.$$ 
Thus, we have 
$$ \#I(j)\le \max I(j) \le (2^{j+1}x_Q)^{1/({1+\tau})} $$ 
if $I(j)$ is nonempty (as otherwise the estimate is trivial) and 
$$\sum_{k\in I(j)} y_k\le (2^{j+1}x_Q)^{\eta/({1+\tau})} \#I(j) \le (2^{j+1}x_Q)^{(1+\eta)/({1+\tau})}.$$ 
From \eqref{eq:auxbx} it follows that 
$$\sum_{k=Q}^R\frac{y_k}{x_k}\le\sum_{j=0}^\infty \frac{(2^{j+1}x_Q)^{(1+\eta)/({1+\tau})}}{2^j x_Q}\ll \frac{1}{x_Q^{(\tau-\eta)/({1+\tau})}}\sum_{j=0}^\infty \biggl(\frac{1}{2^{(\tau-\eta)/({1+\tau})}}\biggr)^j.$$ 
Since $\tau>\eta$, the geometric series converges, giving \eqref{it:Wpoly-pattern}.

Now we turn to Part \eqref{it:Wexp-pattern}.
If $k\in I(j)$, then $e^k\le x_k\le 2^{j+1}x_Q$ and hence, 
$$ k \le \log(2^{j+1}x_Q)\ll j+\log x_Q. $$ 
Therefore, $y_k\le k^\eta\ll (j+\log x_Q)^\eta$ and 
$$\sum_{k\in I(j)} y_k\ll ( j+\log x_Q)^\eta \#I(j)\ll (j+\log x_Q)^{1+\eta}.$$ 
Since $\sum_{j=0}^\infty (j+L)^{1+\eta}/2^j\ll L^{1+\eta}$ for $L\ge 1$, \eqref{eq:auxbx} gives
\[ \sum_{k=Q}^R\frac{y_k}{x_k}\ll\sum_{j=0}^\infty \frac{(j+\log x_Q)^{1+\eta}}{2^j x_Q}\ll\frac{(\log x_Q)^{1+\eta}}{x_Q}.\qedhere \]
\end{proof}

For the rest of this section, we take sequences $\{a_n\}_{n=1}^\infty$ and $\{b_n\}_{n=1}^\infty$, and assume the hypotheses of Theorem~\ref{thm:erdos1051-gen}.
We set 
$c \coloneq c_{\mathbf{w}}$ and
\begin{equation}\label{eq:defofmun}
\mu_n \coloneq  \frac{\log a_n}{c^n},
\end{equation}
so that we have 
\begin{equation}\label{eq:thelimsupcond}
\limsup_{n\to\infty} \mu_n=\infty.
\end{equation}
It will be convenient to shift the indices from the original theorem formulation and denote
\begin{align}
x_n & \coloneq a_{n-d+1}^{w_{0}}a_{n-d+2}^{w_{1}}\cdots a_{n}^{w_{d-1}}=\exp\Bigl(\sum_{j=0}^{d-1} w_j c^{n-d+1+j} \mu_{n-d+1+j}\Bigr), \label{eq:defofxn} \\
y_n & \coloneq b_{n-d+1} \nonumber
\end{align}
for every $n\geq d$, so that
\[ S \coloneq S_{\mathbf{w}}(\{a_n\}_{n=1}^\infty,\{b_n\}_{n=1}^\infty) = \sum_{n=d}^\infty \frac{y_n}{x_n}. \]
From the theorem hypotheses we have
\begin{equation}\label{eq:assumption_on_yn}
y_n \leq n^\eta
\end{equation}
for all $n\geq d$.
We also know that $x_n\geq (n-d+1)^{1+\tau}\geq n^{1+\tau'}$ for $\tau'\in(\eta,\tau)$ and large $n$. By writing $\tau$ in place of $\tau'$ we can immediately assume that
\begin{equation}\label{eq:assumption_on_xn}
x_n \geq  n^{1+\tau}
\end{equation}
for all sufficiently large $n$.
Setting 
\begin{equation}\label{eq:defofdeltan}
\delta_n\coloneq \frac{1}{n^2},
\end{equation}
we see that the peak-index set $\mathcal{P}$ in Lemma~\ref{lem:borel} is infinite and \eqref{eq:thelimsupcond} implies
\begin{equation}\label{eq:thelimsupcond2}
\lim_{\mathcal{P}\ni n\to\infty} \mu_{n+1} = \infty.
\end{equation}

Let us make the following observation regarding ``local peaks''; this is one of the crucial ingredients in our adaptation of the proof from~\cite{Erdos1975}.
In particular, the validity of this inequality heavily depends on our choice of $c$.

\begin{lemma}\label{lem:local-peak}
Suppose we have indices $P, Q$ such that $Q\ge P+d-2$ and 
$$
\mu_{Q+1} > (1+\delta_{Q})\max_{P\leq k\leq Q} \mu_{k}.
$$
Then we have
$$ \frac{(a_P a_{P+1}\cdots a_Q)^W}{x_{Q+1}}\le \Bigl(\frac1{a_{Q+1}}\Bigr)^{1/(Q^2+1)}.
$$
\end{lemma}
\begin{proof}
Let us set $Q_0\coloneq Q-d+2\ge P$. 
Since $c$ is a root of $P_{\mathbf{w}}$, we have
\[ W c^{d-1} = (c-1)\sum_{j=0}^{d-1} w_j c^j, \]
so that
\begin{equation}\label{eq:justauxW}
W\sum_{k=P}^{Q}c^k 
\le W\sum_{k=1}^{Q}c^k
\le \frac{W c^{Q+1}}{c-1}=c^{Q_0}\frac{W c^{d-1}}{c-1}=c^{Q_0}\sum_{j=0}^{d-1} w_j c^{j}.
\end{equation}
Using \eqref{eq:defofmun}, \eqref{eq:defofxn}, \eqref{eq:justauxW}, and \eqref{eq:defofdeltan} we can write
$$
\begin{aligned}
\log \frac{(a_P a_{P+1}\cdots a_Q)^W}{x_{Q+1}}
& = W \sum_{k=P}^Q c^k \mu_k - \sum_{k=Q_0}^{Q+1} w_{k-Q_0} c^k \mu_k\\
& = \sum_{k=P}^{Q_0-1}Wc^k\mu_{k}
+\sum_{k=Q_0}^Q (\underbrace{W-w_{k-Q_0}}_{\geq0})
c^k\mu_{k}
-w_{d-1}c^{Q+1}\mu_{Q+1}\\
&\le\frac{\mu_{Q+1}}{1+\delta_{Q}} 
\Bigl( \underbrace{W\sum_{k=P}^{Q}c^k-c^{Q_0}\sum_{j=0}^{d-2}w_j c^{j}}_{\le w_{d-1}c^{Q+1}}
\Bigr) -w_{d-1}c^{Q+1}\mu_{Q+1}\\
&\le\frac{-\delta_Q\mu_{Q+1}}{1+\delta_{Q}} 
w_{d-1}c^{Q+1}  
\le -\frac{c^{Q+1}\mu_{Q+1}}{Q^2+1}.
\end{aligned}
$$
Exponentiating the last inequality gives $$\frac{(a_P\cdots a_Q)^W}{x_{Q+1}}
\le \exp\Bigl(-\frac{c^{Q+1}\mu_{Q+1}}{Q^2+1}\Bigr)= \exp\Bigl(-\frac{\log a_{Q+1}}{Q^2+1}\Bigr)=a_{Q+1}^{-1/(Q^2+1)},$$ as claimed.
\end{proof}

We also let
$$D_N\coloneq\prod_{k=1}^{N} a_k^W= \exp\Bigl(W \sum_{k=1}^{N} c^k \mu_k\Bigr)$$
and consider the tail sum 
\begin{equation}\label{eq:the_xy_tail}
r_N\coloneq  \sum_{n=N+1}^\infty \frac{y_n}{x_n}. 
\end{equation}
For the purpose of proving Theorem~\ref{thm:erdos1051-gen}, it suffices to establish the following analytic fact; we note that its proof will not require the assumption that $a_i$ or $w_i$ are integers.

\begin{proposition}\label{prop:erdos1051}
Under hypotheses \eqref{eq:assumption_on_yn}, \eqref{eq:assumption_on_xn}, and \eqref{eq:thelimsupcond2}, we have 
\[ \liminf_{N\to\infty} D_N r_N=0. \]
\end{proposition}

\begin{proof}
Let us set
$$ \alpha\coloneq \frac{\tau-\eta}{1+\tau}>0,\quad M\coloneq \Bigl\lceil \frac{2}{\alpha}\Bigr\rceil, $$
and
$$ A\coloneq\max\Bigl\{MW+1,\,\prod_{n=1}^\infty (1+\delta_n),\, W\sup_{x\ge 1}\frac{x^2}{c^x}\Bigr\}. $$
We consider three cases below. The main body of the proof is Case (C), while the other two cases require certain variations of the strategy from~\cite{Erdos1975}. 

\medskip
\noindent\textbf{Case (A).} There exist infinitely many indices $N$ such that $x_{N+1}\ge D_N^M$.
Due to \eqref{eq:assumption_on_yn} and \eqref{eq:assumption_on_xn}, Lemma~\ref{lem:integral-weighted-pattern} Part \eqref{it:Wpoly-pattern} can be applied to \eqref{eq:the_xy_tail} with $Q=N+1$ and letting $R\to\infty$ to give
\[ r_N\ll \frac{1}{x_{N+1}^\alpha}. \]
Hence, for every $N\in\mathbb{N}$ such that $x_{N+1}\ge D_N^M$, 
\[ D_N r_N \ll \frac{D_N}{x_{N+1}^\alpha}\le \frac{x_{N+1}^{1/M}}{x_{N+1}^\alpha}= x_{N+1}^{\,1/M-\alpha}\to 0 \]
as $N\to\infty$, since $1/M\le \alpha/2$ and $x_{N+1}\to\infty$. This proves $\liminf_{N\to\infty} D_N r_N=0$ in Case (A).

\medskip
\noindent\textbf{Case (B).} We have $a_n\ge e^n$ for all sufficiently large $n$.
Take a large $N\in \mathcal{P}$. 
Since 
$$x_n\ge a_n^{w_{d-1}}\ge a_n\ge e^n$$
for all large $n$, we can apply Lemma~\ref{lem:integral-weighted-pattern} Part \eqref{it:Wexp-pattern} with $Q=N+1$ and let $R\to\infty$ to obtain 
$$r_N\ll \frac{(\log x_{N+1})^{1+\eta}}{x_{N+1}}.$$ Also, applying Lemma~\ref{lem:local-peak} with $P=1$ and $Q=N$ yields 
$$\frac{D_N}{x_{N+1}}=\frac{(a_1\cdots a_N)^W}{x_{N+1}}\le \Bigl(\frac1{a_{N+1}}\Bigr)^{1/(N^2+1)}.$$
Combining the last two displays we get
$$\log(D_N r_N)\le-\frac{\log a_{N+1}}{N^2+1}+(1+\eta)\log\log x_{N+1}+O(1).$$ 
Since $x_{N+1}\le a_{N+1}^{\sum_{i=0}^{d-1} w_i}$ and $\log a_{N+1}=c^{N+1}\mu_{N+1}$, we have 
\[ \log\log x_{N+1}\leq \log\log a_{N+1} + O(1)
\ll N+\log\mu_{N+1}. \]
Hence, 
$$\log(D_N r_N)\le-\frac{c^{N+1}\mu_{N+1}}{N^2+1}+O(N)+O(\log\mu_{N+1}),$$ 
which tends to $-\infty$ as $N\to\infty$ along $\mathcal{P}$, because $c>1$ and $\mu_{N+1}\to\infty$. Thus, $D_N r_N\to 0$ along $N\in\mathcal P$, proving the desired result in Case (B).

\medskip
\noindent\textbf{Case (C).} Neither (A) nor (B) holds.
This means
\[ x_{n+1}<D_n^M\ \,\text{for all sufficiently large }n \]
and
\[ a_n<e^n\ \,\text{for infinitely many }n. \] 
We discard those finitely many indices $n$ for which the first of the two estimates fails, i.e., we only work with sufficiently large $n\in\mathbb{N}$.

Because of $a_{n+1}\le x_{n+1}$, we have
$$D_{n+1}=a_{n+1}^W D_n \le (D_n^M)^W D_n = D_n^{MW+1}.$$
Since $A\geq MW+1$, we get by iteration \begin{equation}\label{eq:dnk-pattern}
D_{n+k}\le D_n^{A^k}
\end{equation}
for sufficiently large $n$ and for $k\geq 1$.

Take any $R\in\mathcal{P}$ sufficiently large such that $c^{R+1}>R+1$ and $\mu_{R+1}>1$. 
This implies $\log a_{R+1}=c^{R+1}\mu_{R+1}>R+1$, so $a_{R+1}>e^{R+1}$. Let $P=P(R)<R$ be the largest index with $a_P<e^P$. 
Then 
\begin{equation}\label{eq:a_n_in_between}
a_n\ge e^n\ \,\text{for all } n\in[P+1,R+1].
\end{equation} 
Also, by the standing assumptions,
\[ \lim_{\mathcal P\ni R\to\infty} P(R) =\infty. \]
Moreover,
\begin{equation}\label{eq:DP-pattern}
D_P\le a_P^{PW}\le e^{P^2W}.
\end{equation}
Combined with \eqref{eq:dnk-pattern}, this gives
$$D_{P+k}\le D_P^{A^k}\le e^{P^2W A^k}$$
for all $k\ge0$. Thus, 
$$\mu_{P+k}=\frac{\log a_{P+k}}{c^{P+k}}\le \frac{\log D_{P+k}}{c^{P+k}}\le \frac{P^2W}{c^P}\Bigl(\frac{A}{c}\Bigr)^k.$$ 
Since $A\geq W\sup_{x\ge 1}x^2/c^x$, every large $P$ satisfies
\begin{equation}\label{eq:mupk-pattern}
\mu_{P+k}\le A^{k+1}
\end{equation}
for $k\geq 1$.
In particular, $\mu_{R+1}\le A^{R-P+2}$, so the fact that $\mu_{R+1}\to\infty$ as $R\to\infty$ along $\mathcal P$ implies 
\[ \lim_{\mathcal P\ni R\to\infty}(R-P(R))=\infty. \]
Hence, we may assume $R-P>d$. 

Let us consider the set of ``local peaks,''
$$\mathcal P^\ast\coloneq\Bigl\{m\in[P+d,R]:\ \mu_{m+1}>1\ \text{ and }\mu_{m+1}>(1+\delta_m)\max_{P\le k\le m}\mu_k\Bigr\}.$$ 
Since $R\in\mathcal P$ and $\mu_{R+1}>1$, we have $R\in\mathcal P^\ast$, so $\mathcal{P}^\ast\neq\emptyset$. Let $Q=Q(R)$ be the smallest element of $\mathcal{P}^\ast$. Set temporarily $\bar{\mu}\coloneq \max_{P\le k\le P+d}\mu_k$.
By \eqref{eq:mupk-pattern}, $\bar{\mu}\le A^{d+1}$. By the minimality of $Q$ and the definition of $\mathcal P^\ast$, we inductively obtain 
$$\mu_{P+d+\ell}\le \max\Bigl\{1,\bar{\mu}\prod_{j=P+d}^{P+d+\ell-1}(1+\delta_j)\Bigr\}\le A^{d+1}\prod_{n=1}^\infty(1+\delta_n)\leq A^{d+2}$$ 
for $1\le \ell\le Q-P-d$.
Let us summarize our conclusions as
\begin{equation}\label{eq:mu-max-pattern}
\max_{P\le k\le Q}\mu_k\le A^{d+2},
\quad
\mu_{Q+1}>(1+\delta_Q)\max_{P\le k\le Q}\mu_k,
\quad
\mu_{Q+1}>1,
\end{equation}
and
\[ \lim_{\mathcal P\ni R\to\infty} Q(R) =\infty. \]

Now split 
\begin{equation}\label{eq:thesplitting}
D_Q r_Q
= \underbrace{D_Q\sum_{k=Q+1}^R \frac{ y_k}{x_k}}_{\eqcolon\Sigma_1(R)} + \underbrace{D_Q\sum_{k=R+1}^\infty \frac{ y_k}{x_k}}_{\eqcolon\Sigma_2(R)}.
\end{equation}

\medskip
\noindent\emph{Estimate for $\Sigma_1(R)$.}
Because of \eqref{eq:a_n_in_between}, Lemma~\ref{lem:integral-weighted-pattern} Part \eqref{it:Wexp-pattern} applies and gives
$$\Sigma_1(R)=D_Q\sum_{k=Q+1}^R \frac{y_k}{x_k}\ll \frac{D_Q(\log x_{Q+1})^{1+\eta}}{x_{Q+1}}.$$ 
Thanks to \eqref{eq:mu-max-pattern} we can use Lemma~\ref{lem:local-peak} and combine it with \eqref{eq:DP-pattern} to conclude
$$
\log \frac{D_Q}{x_{Q+1}}
=\log \frac{(a_{P+1}\cdots a_Q)^W}{x_{Q+1}} +\log D_P
\le  -\frac{c^{Q+1}\mu_{Q+1}}{Q^2+1}+P^2W.
$$
By \eqref{eq:mupk-pattern} with $k=Q+1-P$, we have $\mu_{Q+1}\le A^{Q-P+2}$. Hence, $$\log\log x_{Q+1}\le \log(dW\log a_{Q+1})=\log(dW c^{Q+1}\mu_{Q+1}) = O(Q).$$ 
It follows that
$$
\begin{aligned}
\log\Sigma_1(R)
&\le \log \frac{D_Q}{x_{Q+1}} + O(\log\log x_{Q+1}) + O(1)\\
&\le -\frac{c^{Q+1}\mu_{Q+1}}{Q^2+1}+ P^2W+ O(Q)\\
&\le -\frac{c^{Q+1}}{Q^2+1}
+O(Q^2)
\end{aligned}
$$
for sufficiently large $P<Q\leq R$.
The last expression tends to $-\infty$ as $\mathcal P\ni R\to\infty$, since then also $Q=Q(R)\to\infty$.
Hence, $\Sigma_1(R)\to 0$.

\medskip
\noindent\emph{Estimate for $\Sigma_2(R)$.}
This time we can only use \eqref{eq:assumption_on_xn}, so Lemma~\ref{lem:integral-weighted-pattern} Part \eqref{it:Wpoly-pattern} gives 
$$\sum_{k=R+1}^\infty \frac{y_k}{x_k}\ll \frac{1}{x_{R+1}^\alpha}.$$
Thus,
$$\log\Sigma_2(R)\le \log D_Q -\alpha\log x_{R+1}+O(1).$$ 
Using \eqref{eq:mu-max-pattern}, we may write 
$$ \log\frac{D_Q}{D_P}=W\sum_{k=P+1}^Q c^k\mu_k\ll W A^{d+2}\sum_{k=P+1}^Q c^k = O(c^{R}) $$ 
and so, by \eqref{eq:DP-pattern},
$$\log D_Q=\log \frac{D_Q}{D_P}+\log D_P\le
O(c^{R}) + P^2W \le O(c^{R})$$
for large $P<Q\leq R$. 
Since $\log x_{R+1}\ge \log a_{R+1}=c^{R+1}\mu_{R+1}$, 
we see that 
\[ \log\Sigma_2(R)\le c^{R+1}\bigl(-\alpha\mu_{R+1} + O(1)\bigr), \]
which converges to $-\infty$ as $\mathcal{P}\ni R\to\infty$, since $\mu_{R+1}\to\infty$ along $\mathcal{P}$. Hence, $\Sigma_2(R)\to 0$. 

\medskip
From the splitting \eqref{eq:thesplitting} we conclude $\lim_{\mathcal{P}\ni R\to\infty} D_{Q(R)} r_{Q(R)} = 0$.
This completes the proof of Proposition~\ref{prop:erdos1051}, and hence, also that of Theorem~\ref{thm:erdos1051-gen}.
\end{proof}

\begin{remark}
In the above proof, we can also show
\begin{equation}\label{eq:gap-rq}
\lim_{\mathcal P\ni R\to\infty}(R-Q(R))=\infty,
\end{equation}
even though this was not used. An analogous claim was made without proof in~\cite{Erdos1975}.
Indeed, we first note from \eqref{eq:dnk-pattern} that
$$c^{R+1}\mu_{R+1} =\log a_{R+1}\le \log D_{R+1}\le A^{R+1-Q}\log D_Q.$$
By the definition of $D_N$ we have
$$ \log D_Q =\log D_{P+d-1} + W\sum_{k=P+d}^Q c^k\mu_k,$$
while from \eqref{eq:mu-max-pattern} we see
$$
\sum_{k=P+d}^Q c^k\mu_k\le A^{d+2}\sum_{k=P+d}^Q c^k \le A^{d+2}c^{Q+1}/(c-1).
$$
We have seen above that
$ \log D_Q \ll P^2+c^Q$ 
and so
$$c^{R+1}\mu_{R+1} \ll   A^{R+1-Q} (P^2+c^Q).$$
Since $c^{-(R+1)}(P^2+c^Q)$ is bounded and since $\mu_{R+1}\to\infty$, we obtain \eqref{eq:gap-rq}.
\end{remark}

After Proposition~\ref{prop:erdos1051} it is easy to complete the proof of Theorem~\ref{thm:erdos1051-gen}.

\begin{proof}[Proof of Theorem~\ref{thm:erdos1051-gen}]
Each of the numbers $x_d,\ldots,x_N$ clearly divides $D_N$, so 
$$ D_N\Bigl(\frac{y_d}{x_d}+\cdots+\frac{y_N}{x_N}\Bigr)\in\mathbb{Z}. $$ 
Mahler's criterion (Lemma~\ref{lem:mahler}) can be applied to the series with terms
\[ z_n \coloneq \begin{cases}
0 & \text{for } 1\leq n\leq d-1,\\
y_n/x_n & \text{for } n\geq d
\end{cases} \]
and combined with Proposition~\ref{prop:erdos1051} to conclude that $S$ is an irrational number.
\end{proof}

%%%%%

\section{Proof of Theorem~\ref{thm:erdos1051-construction}}
The construction below is an instance of a general ``interval-filling principle,'' which has also been applied recently to several other irrationality problems \cite{Hancl1991}, \cite[\S5]{TY2002}, \cite[\S3]{HT2004}, \cite{Kov25}, \cite[\S5]{KT2025}, \cite[\S2]{CK25}, \cite[\S2]{DK26}.
Namely, suppose that a convergent series $\sum_n \varphi_n(a_1,a_2,\ldots,a_n)$ depends on successive choices $a_n\in J_n$, where the sets $J_n$ are given and finite; typically they are discrete intervals. After fixing a prefix $a_1,\ldots,a_N$, let 
\[ I_N(a_1,\ldots,a_N) \] 
denote the interval between the smallest and largest values that can still be obtained by completing the series tail. Choosing the next term gives ``child'' intervals 
\[ I_{N+1}(a_1,\ldots,a_N,a_{N+1}). \]
If, for all sufficiently large $N$, these child intervals overlap consecutively as $a_{N+1}$ ranges through $J_{N+1}$, then
\[ I_N(a_1,\ldots,a_N) = \bigcup_{a_{N+1}\in J_{N+1}} I_{N+1}(a_1,\ldots,a_N,a_{N+1}). \]
Iterating this identity shows that the set of all attainable sums is a finite union of closed intervals. Equivalently, every point of one of these intervals can be followed through a nested sequence of child intervals, and the shrinking of the series tails then gives an actual completion realizing that point. Once the series sums fill in a non-degenerate interval, we surely know that they attain rational values. 
The same reasoning can also be useful in higher dimensions \cite[\S7]{KT2025} and even in infinite-dimensional spaces \cite[\S8]{KT2025}.

The following is the main technical lemma, with a non-constructive proof that we just described. Let us fix $d$ and $\mathbf{w}$, and briefly write $S(\{a_n\}_{n=1}^{\infty})$ instead of 
$S_{\mathbf{w}}(\{a_n\}_{n=1}^\infty,\{1\}_{n=1}^\infty)$.

\begin{lemma}\label{lm:negative}
Let $\{\beta_n\}_{n=1}^\infty$ and $\{\gamma_n\}_{n=1}^\infty$ be monotonically increasing sequences of positive integers such that $\beta_n\leq \gamma_n$ for every $n\in\mathbb{N}$,
\begin{equation}\label{eq:sumcd}
\sum_{n=1}^\infty \frac{1}{\beta_n^{w_0} \beta_{n+1}^{w_1} \cdots \beta_{n+d-1}^{w_{d-1}}} < \infty, 
\end{equation}
\begin{equation}\label{eq:limcd}
\lim_{n\to\infty} \frac{\beta_n}{\gamma_n} = 0, 
\end{equation}
and
\begin{equation}\label{eq:cond_limit}
\lim_{n\to\infty} \Bigl(\frac{\beta_{n}}{\gamma_{n+1}}\Bigr)^{w_0} \cdots \Bigl(\frac{\beta_{n+d-2}}{\gamma_{n+d-1}}\Bigr)^{w_{d-2}} \Bigl(\frac{\beta_{n+d-1}}{\beta_{n+d}}\Bigr)^{w_{d-1}} \beta_{n+d-1} = \infty.
\end{equation}
Then the set of infinite sums
\begin{equation}\label{eq:set_sums} 
\bigl\{ S(\{a_n\}_{n=1}^\infty) \,:\, (\forall n\in\mathbb{N})(a_n\in [\beta_n,\gamma_n]\cap\mathbb{N}) \bigr\} 
\end{equation}
is equal to a finite union of non-degenerate closed bounded intervals. In particular, it contains a rational number.
\end{lemma}

\begin{proof}
Let $J_n \coloneq [\beta_n,\gamma_n]\cap\mathbb{N}$ for every $n\in\mathbb{N}$.
We only need to show that there exists an $M\in\mathbb{N}$ such that, for all indices $N\geq M$ and all choices of $a_n\in J_n$, $1\leq n\leq N$, the closed interval
\begin{equation}\label{eq:interval0}
\bigl[ S(a_1,\ldots,a_N,\gamma_{N+1},\gamma_{N+2},\ldots),\, S(a_1,\ldots,a_N,\beta_{N+1},\beta_{N+2},\ldots) \bigr] 
\end{equation}
is equal to the union of
\begin{equation}\label{eq:interval1}
\bigl[ S(a_1,\ldots,a_N,a_{N+1},\gamma_{N+2},\ldots),\, S(a_1,\ldots,a_N,a_{N+1},\beta_{N+2},\ldots) \bigr] 
\end{equation}
taken over $a_{N+1}\in J_{N+1}$.
Afterwards, we can conclude that \eqref{eq:set_sums} is equal to the union of intervals
\begin{equation}\label{eq:interval2}
\bigl[ S(a_1,\ldots,a_M,\gamma_{M+1},\gamma_{M+2},\ldots),\, S(a_1,\ldots,a_M,\beta_{M+1},\beta_{M+2},\ldots) \bigr] 
\end{equation}
over the finitely many choices of $a_n\in J_n$, $1\leq n\leq M$.
Indeed, for an arbitrary real number $x$ from \eqref{eq:interval2} we use the claim to inductively construct $a_n\in J_n$, $n=M+1,M+2,\ldots$, such that $x$ belongs to \eqref{eq:interval0} for all $N\geq M$.
Since $S(\{a_n\}_{n=1}^\infty)$ also lies in the intervals \eqref{eq:interval0} and their lengths shrink to zero as $N\to\infty$ thanks to \eqref{eq:sumcd}, we obtain the representation $x=S(\{a_n\}_{n=1}^\infty)$, proving the lemma.

To complete the proof it remains to show that, for all sufficiently large $N$ and all $a_n\in J_n$, $1\leq n\leq N+1$, $a_{N+1}\neq\gamma_{N+1}$, the left endpoint of \eqref{eq:interval1} does not exceed the right endpoint of \eqref{eq:interval1} with $a_{N+1}$ replaced by $a_{N+1}+1$, i.e.,
\begin{equation}\label{eq:endpoints1}
S(a_1,\ldots,a_N,a_{N+1},\gamma_{N+2},\ldots) \leq S(a_1,\ldots,a_N,a_{N+1}+1,\beta_{N+2},\ldots).
\end{equation}
This will guarantee that the intervals \eqref{eq:interval1} fully cover \eqref{eq:interval0}.
Observe that
\[ \frac{1}{\gamma_n^{w_0} \gamma_{n+1}^{w_1} \cdots \gamma_{n+d-1}^{w_{d-1}}} \leq \frac{1}{\beta_n^{w_0} \beta_{n+1}^{w_1} \cdots \beta_{n+d-1}^{w_{d-1}}} \]
holds for all $n\geq N+2$ and that
\[ \frac{1}{a_n^{w_0} \cdots a_{N+1}^{w_{N-n+1}} \gamma_{N+2}^{w_{N-n+2}} \cdots \gamma_{n+d-1}^{w_{d-1}}} \leq \frac{1}{a_n^{w_0} \cdots (a_{N+1}+1)^{w_{N-n+1}} \beta_{N+2}^{w_{N-n+2}} \cdots \beta_{n+d-1}^{w_{d-1}}} \]
holds for sufficiently large $N$ and for $N-d+4\leq n\leq N+1$, simply due to $w_{d-1}>0$ and $\lim_{n\to\infty}\beta_{n+d-1}/\gamma_{n+d-1}=0$.
Thus, inequality \eqref{eq:endpoints1} will follow from comparing contributions of only two terms from each of the two series, those corresponding to indices $n=N-d+2$ and $n=N-d+3$:
\begin{align*}
& \frac{1}{a_{N-d+2}^{w_0} \cdots a_N^{w_{d-2}} a_{N+1}^{w_{d-1}}} 
+ \frac{1}{a_{N-d+3}^{w_0} \cdots a_N^{w_{d-3}} a_{N+1}^{w_{d-2}} \gamma_{N+2}^{w_{d-1}}} \\
& \leq \frac{1}{a_{N-d+2}^{w_0} \cdots a_N^{w_{d-2}} (a_{N+1}+1)^{w_{d-1}}}
+ \frac{1}{a_{N-d+3}^{w_0} \cdots a_N^{w_{d-3}} (a_{N+1}+1)^{w_{d-2}} \beta_{N+2}^{w_{d-1}}}.
\end{align*}
Rearranging the terms, this inequality can be equivalently rewritten as
\begin{align*}
& \Bigl(\frac{a_{N-d+2}}{a_{N-d+3}}\Bigr)^{w_0} \cdots \Bigl(\frac{a_{N}}{a_{N+1}}\Bigr)^{w_{d-2}} \Bigl(\frac{a_{N+1}}{\beta_{N+2}}\Bigr)^{w_{d-1}} a_{N+1} \\
& \geq a_{N+1} \biggl( 1 - \Bigl(\frac{a_{N+1}}{a_{N+1}+1}\Bigr)^{w_{d-1}} \biggr)
\bigg/ \biggl( \Bigl(\frac{a_{N+1}}{a_{N+1}+1}\Bigr)^{w_{d-2}} - \Bigl(\frac{\beta_{N+2}}{\gamma_{N+2}}\Bigr)^{w_{d-1}} \biggr).
\end{align*}
On the one hand, by \eqref{eq:limcd} and the limit $$\lim_{x\to\infty} x\left(1-\Bigl(\frac{x}{x+1}\Bigr)^{w_{d-1}}\right) = w_{d-1},$$ the right-hand side is bounded from above by a finite constant that depends only on $\mathbf{w}$.
On the other hand, the left-hand side is at least
\[ \Bigl(\frac{\beta_{N-d+2}}{\gamma_{N-d+3}}\Bigr)^{w_0} \cdots \Bigl(\frac{\beta_{N}}{\gamma_{N+1}}\Bigr)^{w_{d-2}} \Bigl(\frac{\beta_{N+1}}{\beta_{N+2}}\Bigr)^{w_{d-1}} \beta_{N+1}, \]
which, by \eqref{eq:cond_limit}, converges to infinity as $N\to\infty$. Thus, the left-hand side is larger than the right-hand side for all sufficiently large $N\in\mathbb{N}$, and all allowed choices of $a_n$, $1\leq n\leq N+1$.
\end{proof}

Lemma~\ref{lm:negative} does not appear in the literature, but, as we have already mentioned, somewhat similar constructions have been used in the context of irrationality of the Ahmes and Cantor series.

\begin{proof}[Proof of Theorem~\ref{thm:erdos1051-construction}]
Denote $c\coloneq\tilde{c}_{\mathbf{w}}$, so that
\begin{equation}\label{eq:eqforc}
(c-1)\sum_{j=0}^{d-1} w_j c^j - c^{d-1} = 0.
\end{equation}
For a given $C>1$, one can choose, say,
\[ \beta_n = \lfloor C^{c^n+n^2+1} \rfloor\quad\text{and}\quad \gamma_n = \lfloor C^{c^n+n^2+n} \rfloor \]
for every $n\in\mathbb{N}$.
Conditions \eqref{eq:sumcd} and \eqref{eq:limcd} are clearly satisfied, while the limit in \eqref{eq:cond_limit} is actually
\[ \lim_{n\to\infty} C^{\sum_{j=0}^{d-1}w_j (c^{n+j}-c^{n+j+1}) + c^{n+d-1} + n^2 + O(n)} = \lim_{n\to\infty} C^{n^2 + O(n)} = \infty, \]
where we used \eqref{eq:eqforc}, i.e., our choice of $c$.

Lemma~\ref{lm:negative} gives a sequence of positive integers $\{a_n\}_{n=1}^\infty$ satisfying $\beta_n \leq a_n\leq \gamma_n$ for every $n\in\mathbb{N}$ and $S(\{a_n\}_{n=1}^\infty)\in\mathbb{Q}$.
Such a sequence clearly grows according to \eqref{eq:examplegrowth}.
A minor detail is that it only becomes strictly increasing from some index onward; we surely have $a_n<a_{n+1}$ for all $n$ large enough that $\gamma_n<\beta_{n+1}$.
However, modifying finitely many terms of $\{a_n\}_{n=1}^\infty$ does not affect the rationality of the sum $S(\{a_n\}_{n=1}^\infty)$. Thus, the double-exponential growth of this sequence leaves more than enough room to redefine its first finite block as $a_n\coloneq n$ and end up with a strictly increasing sequence that still satisfies \eqref{eq:examplegrowth}.
\end{proof}

\section{Remaining open questions}
The positive result in Theorem~\ref{thm:erdos1051-gen} involves the root $c_{\mathbf{w}}$ of
\[ (x-1)\sum_{j=0}^{d-1}w_jx^j - (\max_j w_j) x^{d-1} = 0, \]
whereas the construction in Theorem~\ref{thm:erdos1051-construction} involves the largest positive root $\tilde{c}_{\mathbf{w}}$ of
\[ (x-1)\sum_{j=0}^{d-1}w_jx^j - x^{d-1} = 0. \]
These constants coincide when $w_j\in\{0,1\}$ for all $j$, and in those cases the theorems give a sharp threshold. Otherwise, $c_{\mathbf{w}}$ and $\tilde{c}_{\mathbf{w}}$ may differ, and the authors are currently undecided on which of them is closer to the true threshold for the (ir)rationality of $\sum_n 1/(a_n^{w_0}\cdots a_{n+d-1}^{w_{d-1}})$.

Numerous problems of Erd\H{o}s on the irrationality of Ahmes series are still open. We can briefly mention just a couple of those that are concerned with general series of a certain type, rather than with concrete numerical values.
Erd\H{o}s asked about the fastest growth of a sequence of positive integers $\{a_n\}_{n=1}^\infty$ such that both the series $\sum_n 1/a_n$ and $\sum_n 1/(a_n+1)$ have rational sums \cite[p.~64]{ErdosGraham1980}, \cite[p.~104]{Erd1988}. It is now known that double-exponential growth is possible \cite[Cor.~2.9]{KT2025}, but the precise exponent $c>1$ in the threshold condition $\lim_{n\to\infty}a_n^{1/c^n}=\infty$ remains open.
Erd\H{o}s and Graham also suspected that the series
\[ \sum_{n=1}^{\infty}\frac{1}{(n+1)(n+2)\cdots(n+f(n))} \]
has an irrational sum whenever the number of factors in the denominator, $f(n)$, grows to infinity \cite[p.~65]{ErdosGraham1980}. This is now known to be false in general \cite[Thm.~1]{CK25}, but a sub-question on monotonically increasing functions $f$ remains unanswered.

\section*{Acknowledgements}
The third named author is supported by Mid-Career Researcher Program (RS-2023-00278510) through the National Research Foundation funded by the government of Korea, and by KIAS Individual Grant (MG073602) at Korea Institute for Advanced Study. The third named author is also supported by KIAS--KAIST Joint Research Group Grant.

The fourth named author is supported in part by the Croatian Science Foundation under the project HRZZ-IP-2022-10-5116 (FANAP), and in part by the European Union -- NextGenerationEU through the National Recovery and Resilience Plan 2021--2026, via an institutional grant of the University of Zagreb Faculty of Science IK IA 1.1.3. Impact4Math.

We thank Thomas Bloom for creating and maintaining the website \emph{Erd\H{o}s problems} and Terence Tao for valuable comments on an earlier draft of this paper.
We thank the authors of~\cite{aletheia} for building the \emph{Aletheia} agent and for allowing the authors of this work to utilize it.

\end{document}